# Unique polynomial solution of $m/n = 1/x + 1/y + 1/z$ for $n \equiv b \mod a$ if $(a,m)=1$.

Bernd R. Schuh

Dr. Bernd Schuh, D-50968 Cologne, Germany, bernd.schuh@netcologne.de

**Abstract**

Necessary and sufficient conditions for the existence of an integer solution of the diophantine equation $m/(b+a\lambda) = 1/x(\lambda) + 1/y(\lambda) + 1/z(\lambda)$ are explicitely given for $a,b$ coprime and $a$ not a multiple of $m$. The solution has the form $x(\lambda) = kn(\lambda)$, $z(\lambda) = (kl/r)(s+r\lambda)$, $y(\lambda) = n(\lambda)(s+r\lambda)$ where parameters $l,k,s,r \in \mathbb{Z}_+$ obey certain conditions depending on $a,b$. The conditions imply restrictions for some choices of $a,b$, which differ from the ones known in the case $m=4$. E.g., the modulus must be of the form $l(mk-1)$. One can also deduce, that primes of the form $4K+1$ are excluded as modulus. Also if $a = p \neq m$ is prime and $b = a+1$, i.e. $n \equiv 1 \mod p$, polynomial solutions are shown to be impossible. All results are valid for integers $m \geq 4$.

**I Introduction**

Two well known conjectures by Erdös-Straus and Sierpinski state that the diophantine equation

$$m/n = 1/x + 1/y + 1/z \qquad (1)$$

has integer solutions $x, y, z$ for every integer $n \geq 2$ and $m = 4$ (Erdös-Straus) or $m = 5$ (Sierpinski). There is an impressive body of evidence for the validity of both conjectures, see e.g. [1,2,3,5], but no valid proof.

In the following we will consider the residue class $n \equiv n_0 \mod n_1$ and polynomial denominators in (1):

$$\frac{m}{n_0 + n_1\lambda} = \frac{1}{x(\lambda)} + \frac{1}{y(\lambda)} + \frac{1}{z(\lambda)} \qquad (2)$$

A number of solutions of this equation is known for certain values of $m$, but there is no covering set of equations with $n_1$ fixed and $n_0$ running through all residues modulo $n_1$, see e.g.[4]. To the contrary, for $m \equiv 0 \bmod 4$ Schinzel [6] has shown that solutions of the form (2) are impossible if $n_0$ is a quadratic residue modulo $n_1$. This strong result is limited to the case $m = 4$, however. A simple counterexample for $m = 5$ is $n \equiv 7 \bmod 9$. 7 is a quadratic residue modulo 9, but (2) is solved by $x(\lambda) = 2n;\ z(\lambda) = 2(1+\lambda);\ y(\lambda) = (1+\lambda)n$.

Most investigations in the case $m = 4$ center around a modulus divisible by 4. Thus the main assumption of our investigation is not satisfied in this case. But for $m > 4$ a fresh look might be promising. We will see, however, that severe restrictions for the existence of integer polynomial solutions to equation (2) arise in case $m$ does not divide $n_1$. We derive such restrictions from a unique parametric solution to equation (2). The main result can be cast into the following theorem.

**Theorem**. Let $m, n_0, n_1$ be positive integers with $m \geq 4$, $(n_0, n_1) = 1$ and $(n_1, m) = 1$. Equation (2) has a unique solution with integer polynomials $x(\lambda), y(\lambda), z(\lambda)$ given by

$$x(\lambda) = kn(\lambda)$$
$$z(\lambda) = (kl/r)(s + r\lambda) \qquad (3)$$
$$y(\lambda) = n(\lambda)(s + r\lambda)$$

if and only if there exist parameters $l, k, s, r \in \mathbb{Z}_+$ such that conditions i) – iii) are satisfied:

i) $\qquad n_1 = l(mk - 1)$

ii) $\qquad sn_1 = kl + rn_0$

iii) $\qquad skl/r \in \mathbb{Z}_+$

As will become clear in the course of the proof, these conditions phrased in terms of the parameters $x_0, y_0, z_0$ solving equation (2) with $\lambda = 0$ simply state, among other things, that two of these parameters, say $x_0, y_0$ have to be multiples of $n_0$, $x_0 = kn_0, y_0 = sn_0$ and fulfill $y_0((mz_0 - n_0)x_0 - n_0 z_0) = n_0 x_0 z_0$ which is an equivalent form of equation (2) and condition ii) with appropriate definitions.

Before we prove the theorem let us collect some consequences. An obvious one is the following.

**Corollary 1**. An integer solution of equation (2) is not compatible with $(n_0, n_1) = 1$ and $n_1$ being <u>not</u> a multiple of $4k - 1$ for some $k \in \mathbb{N}$ <u>and not</u> a multiple of $m$.

Another obvious consequence is

**Corollary 2**. There are no integer polynomial solutions to equation (2) if $n_1$ is prime and $(n_0, n_1) = 1$, unless $n_1 = m$ or $n_1 = 4k - 1$ for some $k \in \mathbb{N}$.

In particular

**Corollary 3**. Integer polynomial solutions of (2) do not allow primes of the form $4K+1$ as modulus if $(n_1, m) = 1$.

**Corollary 4**. If $n \equiv 1 \bmod p$ where $p$ is a prime and $(p, m) = 1$ there is no integer polynomial solution to (2).

Proof. It is clearly equivalent to prove the corollary for $n \equiv (p+1) \bmod p$. The main assumptions of the theorem are valid, i.e. neither $n_1$ and $n_0$ nor $n_1$ and $m$ have a common divisor except 1. We assume, an integer polynomial solution of equation (2) exists. Then, according to the theorem parameters $r, s, k, l$ exist such that $p = n_1 = l(mk-1)$. Here $l = 1$ since $p$ is prime. Then from condition ii) one gets $s(mk-1) = k(mr+1)$. Since $(k, mk-1) = 1$, $s/k = (mr+1)/(mk-1)$ must be a positive integer. There are two integer solutions to this equation:

    a.     $s/k = m-1$ with $k = 1+\lambda$ and $r = m-2+(m-1)\lambda$

    b.     $s/k = m-1+m\lambda$ with $k = 1$ and $r = m-2+(m-1)\lambda$

with $\lambda \in \mathbb{N}$. In case a. one has
$sk/r = k^2(m-1)/(k(m-1)-1) = k(1+1/(k(m-1)-1)) = (1+\lambda) + (1+\lambda)/(m-2+(m-1)\lambda)$. In case b.
$sk/r = (m-1+m\lambda)/(m-2+(m-1)\lambda) = 1+(1+\lambda)/(m-2+(m-1)\lambda)$. Thus in both cases $skl/r \notin \mathbb{N}$, since $(1+\lambda) < (m-2)(1+\lambda) \leq m-2+(m-1)\lambda$, a violation of condition iii). Since the necessary set of integers to fulfill (2) with positive integer coefficients does not exist, such a solution does not exist in this case.

**II Proof of the theorem.**

To show that the existence of parameters $s, k, l, r$ with properties i) – iii) is sufficient for the existence of an integer polynomial solution of equation (2) one defines $x(\lambda), y(\lambda), z(\lambda)$ via equation (3). Because of condition iii) all coefficients are positive integers. Then one calculates, using i) and ii)

$$1/x + 1/y + 1/z$$
$$= \frac{l(s+r\lambda) + rn_0 + rn_1\lambda + kl}{kl(s+r\lambda)n}$$
$$= \frac{ls + lr\lambda + sn_1 + rn_1\lambda}{kl(s+r\lambda)n} \quad \text{(use ii) to get)}$$

$$= \frac{(s+r\lambda)(l+n_1)}{kl(s+r\lambda)n} \quad \text{(use i) to get)}$$
$$= \frac{lkm}{lkn} = m/n$$

as required. The solution (3) is by no means a lucky choice, but up to possibly different choices of the parameters $s, k, l, r$ unique indeed, as we will see in proving the necessity of the conditions i) – iii).

We do so by first proving that two of the three polynomials $x, y, z$ necessarily are linear in the variable $\lambda$ whereas the third has degree 2. This proof is done via 5 lemmas.

**Lemma 1**. If (2) has a solution in $\mathbb{Q}[\lambda]$ then at least one of the polynomials $x$, $y$, or $z$ has degree 1.

Proof. Rewrite (2) as $mxyz = n\{xy + xz + yz\}$ ($\lambda$-dependence surpressed) and let w.l.o.g. the degrees of the polynomials be ordered as $d(z) \leq d(x) \leq d(y)$ then the degree of the polynomial on the l.h.s. is $d(l.h.s.) = d(z) + d(x) + d(y)$ whereas on the r.h.s. one has $d(r.h.s.) = 1 + d(x) + d(y)$. Equating the two degrees confirms the statement.

**Lemma 2**. If two polynomials in (2) have degree 1, then the degree of the third polynomial cannot be larger than 3.

Proof. Let w.l.o.g. $d(z) = d(x) = 1$ and rewrite (2) as

$$y[mxz - n(x+z)] = nxz. \qquad (2')$$

Then obviously $d(y) \leq d(l.h.s.) = d(r.h.s.) = 1 + 1 + 1 = 3$.

**Lemma 3**. Let (2) be valid. Set $z(\lambda) = z_0 + z_1\lambda$. If both polynomials $x(\lambda), y(\lambda)$ have degrees larger than 1 then

$$z_1 = n_1 / m \qquad (4)$$

Proof. According to lemma 1 one of the three polynomials on the r.h.s. of (2) has degree 1. Let this polynomial be $z(\lambda)$. Now consider (2) for $\lambda \to \infty$ and compare the leading terms on both sides of the equation: $\frac{m}{n_1\lambda}(1 + O(\lambda^{-1})) = \frac{1}{z_1\lambda}(1 + O(\lambda^{-1})) + O(\lambda^{-\min\{d(x), d(y)\}})$. Since $\min\{d(x), d(y)\} \geq 2$ by assumption, equation (3) follows.

Suppose $x(\lambda), y(\lambda), z(\lambda)$ are positive integer solutions of (2). Assume first that two of the polynomials have degrees larger than 1. Then the third must have degree 1 according to lemma 1. But according to lemma 3 the leading coefficient of this polynomial fulfills (4). Since $z_1$ is assumed to be integer we have $n_1 \equiv 0 \bmod m$. But that is contrary to the main assumption of the theorem $(m, n_1) = 1$. So we can dismiss the possibility that two polynomials have degrees larger than 1. Now there are only three possibilities left for the degrees of integer polynomials solving equation (2):

A. $d(x) = d(y) = d(z) = 1$ ; B. $d(x) = d(z) = 1, d(y) = 2$ ; C: $d(x) = d(z) = 1, d(y) = 3$

In the next step we will show that possibilities A and C can be discarded.

Let us consider these cases separately.

**Lemma 4**. Let $n(\lambda) = n_0 + n_1\lambda$ and suppose that all three polynomials $x, y, z$ have degree 1. If we write $x(\lambda) = x_0 + x_1\lambda$, $y(\lambda) = y_0 + y_1\lambda$, $z(\lambda) = z_0 + z_1\lambda$ then the only solution of (2) reads

$$x_1 = x_0 n_1 / n_0; \quad y_1 = y_0 n_1 / n_0; \quad z_1 = z_0 n_1 / n_0.$$

Proof. Multiply equ. (2) by $\lambda$ and write l.h.s. and r.h.s. as an expansion in $1/\lambda$. The result is

$$\frac{m\lambda}{n} = \frac{m}{n_1}\{1 - \frac{n_0}{n_1}\frac{1}{\lambda} + (\frac{n_0}{n_1})^2 \frac{1}{\lambda^2} - ...\} = \frac{1}{x_1}\{1 - \frac{x_0}{x_1}\frac{1}{\lambda} + (\frac{x_0}{x_1})^2 \frac{1}{\lambda^2} - ...\} + (x \to y) + (x \to z)$$

Comparing coefficients of the same powers in $1/\lambda$ yields

$$(\frac{n_0}{n_1})^r = \frac{n_1/m}{x_1}(\frac{x_0}{x_1})^r + \frac{n_1/m}{y_1}(\frac{y_0}{y_1})^r + \frac{n_1/m}{z_1}(\frac{z_0}{z_1})^r$$

valid for every $r \in \mathbb{N}'$. Subtracting the zeroth order term $1 = \frac{n_1}{m}(\frac{1}{x_1} + \frac{1}{y_1} + \frac{1}{z_1})$ one gets

$$0 = (1 - (\frac{x_0 n_1}{n_0 x_1})^r)\frac{1}{x_1} + (1 - (\frac{y_0 n_1}{n_0 y_1})^r)\frac{1}{y_1} + (1 - (\frac{z_0 n_1}{n_0 z_1})^r)\frac{1}{z_1}$$

For all $r \in \mathbb{N}'$. None of the terms $\frac{x_0 n_1}{n_0 x_1}, \frac{y_0 n_1}{n_0 y_1}, \frac{z_0 n_1}{n_0 z_1}$ can be smaller or larger than zero, because then the equation cannot be fulfilled for all $r$. Thus $1 = \frac{x_0 n_1}{n_0 x_1} = \frac{y_0 n_1}{n_0 y_1} = \frac{z_0 n_1}{n_0 z_1}$ as stated.

In order to yield an integer solution of this kind and also to fulfill $(n_0, n_1) = 1$, $n_0$ must divide all three $x_0, y_0, z_0$. But then one gets for the zeroth order solution of (2):

$m = n_0/x_0 + n_0/y_0 + n_0/z_0 = 1/a + 1/b + 1/c$ with positive integers $a, b, c$, a contradiction to $m \geq 4$. So via Lemma 4 we can dismiss case A.

**Lemma 5**. There is no integer (even no real) polynomial solution of (2) if two polynomials have degree 1 and the third has degree 3 and $x_0, y_0, z_0$ are positive, rational numbers that solve equation (2) for $\lambda = 0$.

Proof. Let w.l.o.g. $y(\lambda) = y_0 + y_1\lambda + y_2\lambda^2 + y_3\lambda^3$ and $x(\lambda) = x_0 + x_1\lambda$, $z(\lambda) = z_0 + z_1\lambda$. Write (2) in the form $xy(mz - n) = nz(x + y)$ ($\lambda$-dependence surpressed) and compare coefficients of equal powers of $\lambda$. In fourth and fifth order this yields the two equations

$$(mz_1 - n_1)x_1 = n_1 z_1$$

$$(m_0 z_0 - n_0)x_1 + (mz_1 - n_1)x_0 = n_0 z_1 + n_1 z_0$$

These two equations determine $x_1$ and $z_1$. From these, introducing abbreviations

$$\chi := mx_1 - n_1, \quad \zeta := mz_1 - n_1 \quad \overline{z_0} := mz_0 - n_0, \quad \overline{x_0} := mx_0 - n_0$$

one gets the two equations

$$\chi\zeta = n_1^2$$
$$\overline{z_0}\chi + \overline{x_0}\zeta = 2n_0 n_1$$

These can be solved for $\chi$ with solution

$$\overline{z_0}\chi = n_0 n_1 \pm n_1 \sqrt{n_0^2 - \overline{x_0 z_0}}$$

Now $n_0^2 - \overline{x_0}\overline{z_0} = mn_0(x_0 + z_0) - m^2 x_0 z_0 = mn_0 x_0 z_0 \left[1/z_0 + 1/x_0 - \dfrac{m}{n_0}\right] = -mn_0 x_0 z_0 / y_0 < 0$. Thus there is no real-valued solution to $x_1$ and $z_1$ and no real-valued polynomial solution if the degrees are as assumed in the lemma.

The only possibility left for the form of the solution is the one reflected in solution (3): two polynomials have degree 1, the third has degree 2. Assuming this we can now show

**Lemma 6**. There are two distinct rational polynomial solutions of (2) if two polynomials have degree 1 and the third has degree 2 and $x_0, y_0, z_0$ are positive, rational numbers that solve equation (2) for $\lambda = 0$. These solutions, in an obvious notation, are given by:

$$z_+(\lambda) = z_0 + \dfrac{n_1}{n_0}\dfrac{y_0 z_0}{y_0 + z_0}\lambda \qquad (5a)$$

$$x_+(\lambda) = x_0 + x_0 \dfrac{n_1}{n_0}\lambda \qquad (5b)$$

$$y_+(\lambda) = y_0 + y_0 \dfrac{n_1}{n_0}(1 + \dfrac{y_0}{y_0 + z_0})\lambda + (y_0 \dfrac{n_1}{n_0})^2 \dfrac{1}{y_0 + z_0}\lambda^2 \qquad (5c)$$

and

$$z_-(\lambda) = z_0 + \dfrac{n_1}{m}\dfrac{x_0 + z_0}{x_0}\lambda \qquad (6a)$$

$$x_-(\lambda) = x_0 + \dfrac{n_1}{m}\dfrac{x_0 + z_0}{z_0}\lambda \qquad (6b)$$

$$y_-(\lambda) = y_0 + (2\dfrac{y_0 n_1}{n_0} - \dfrac{n_1}{m})\lambda + \dfrac{n_1}{n_0}(\dfrac{y_0 n_1}{n_0} - \dfrac{n_1}{m})\lambda^2 \qquad (6c).$$

Proof. We choose $y$ to be the polynomial with degree $d = 2$. According to equation (2') the polynomial $X = mzx - n(x+z)$ must have degree 1. Thus $X_2 = 0$. So equation (2') can be written

$$y/y_0 = \dfrac{(1 + n_1'\lambda)(1 + x_1'\lambda)(1 + z_1'\lambda)}{(1 + X_1'\lambda)} \qquad (7a)$$

where $X_0 = n_0 x_0 z_0 / y_0$ has been used and all polynomials appear normalized to their constant value, $X' = X/X_0$, $n_1' = n_1/n_0$ etc.. Equation (7a) has two distinct solutions, namely

$$X_1/X_0 = n_1/n_0 \qquad (7b)$$

$$X_1/X_0 = z_1/z_0 \qquad (7c)$$

The third possibility $X_1' = x_1'$ will not lead to a distinct solution since it amounts to an interchange of variables $x$ and $z$, only.

Finally one needs to determine $z_1$ and $x_1$ for each of the choices (7b), (7c) from the two equations

$$X_1 = (mx_0 - n_0)z_1 + (mz_0 - n_0)x_1 - n_1(x_0 + z_0) \quad (8a)$$

$$X_2 = mx_1 z_1 - n_1(x_1 + z_1) = 0 \quad (8b)$$

We outline the calculation for the case (7c). With abbreviations $\sigma = x_0 + z_0$, $p = x_1 n_0 / n_1$, $q = z_1 n_0 / n_1$ and $\alpha = m/n_0$ equations (8a,b) read:

$$\alpha pq = p + q \quad \text{and} \quad \sigma = qx_0/z_0 + p(\alpha z_0 - 1)$$

These can be transformed into a quadratic equation for either variable. Its two solutions are:

$$q_+ = z_0 \quad q_- = \sigma/(\alpha x_0) \quad \text{and} \quad p_+ = x_0 y_0/(x_0 + y_0) \quad p_- = \sigma/(\alpha z_0).$$

Multiplication with $n_1/n_0$ yields the first order coefficients of $x(\lambda), z(\lambda)$ as stated in equations (5) and (6).

Finally, $y_1$ and $y_2$ are determined from equation (7a) with these solutions inserted. The result is as stated in equs. (5) and (6). A similar calculation for $X_1' = n_1/n_0$ does not lead to different solutions. This completes the proof of lemma 6.

More details can be found in [7], where a complete solution of (2) in polynomials with positive rational coefficients and no restrictions on $n_1, n_0$ is given.

So far we have shown: if an integer solution to (2) exists and also $(n_0, n_1) = 1$ and $(n_1, m) = 1$ are valid this solution can only have one of the forms given in (5) or (6). Now we can turn to the last step in proving the theorem by deducing conditions i) – iii). Consider the "-"-solution first. Let $y_1$ and $y_2$ be the first und second order coefficients of $y_-$. From (6c) one calculates $n_1 y_1 - n_0 y_2 = n_1^2 y_0 / n_0$. For an integer solution the l.h.s. must be an integer. Thus $y_0/n_0 \in \mathbb{Z}$ since $(n_0, n_1) = 1$. Thus also $2n_1 y_0 / n_0 - y_1 = n_1/m$ should be an integer. But that violates the assumption of the theorem. Therefore this solution is not an option.

Consider the "+"-solution next, equation (5). Since $(n_0, n_1) = 1$ the first order coefficient $x_1$ in (5b) is a positive integer if and only if $x_0/n_0 = k \in \mathbb{N}$. Then from (5a), using $m/n_0 = 1/x_0 + 1/y_0 + 1/z_0$:

$z_1 = n_1/(m - 1/k)$ or $mz_1 = mkn_1/(mk-1) = n_1 + n_1/(mk-1)$. If $n_1 < mk - 1$, $mz_1$ and thus $z_1$ cannot be an integer as is assumed in the theorem. Thus the only possibility left is $n_1 = l(mk-1)$ with positive integers $k, l$, as stated in condition i) of the theorem. Next calculate from (5c) $y_1 n_1 - y_2 n_0 = n_1^2 y_0 / n_0$. This can be an integer if and only if $y_0 = s n_0$ with a positive integer $s$. Now

$y_1 = (sn_1^2 + n_0 y_2)/n_1 = sn_1 + n_0 lks/z_0$ and $y_2 = n_1 lks/z_0$. Therefore $y_1$ and $y_2$ can be positive integers if and only if $lks/z_0$ is an integer; we call it $r$. So we have $y_1 = sn_1 + rn_0$. On the other hand one can calculate $y_1 = 2sn_1 - sn_1 z_0 / (y_0 + z_0) = 2sn_1 - slx_0 / y_0 = 2sn_1 - kl$. These two equations for $y_1$ establish condition ii). Thus we have derived conditions i) to iii) from the existence of an integer polynomial solution of equation (2). Expressing the coefficients in equation (5) by the newly introduced parameters one gets exactly the solution (3). Thus this solution is unique under the given premises up to possibly different choices of the parameters.

**III Summary**.

One main point of this paper is: If $n_1$ is not a multiple of $m$, then a necessary condition for the existence of integer polynomial solutions to equation (2) is that $n_1 = l(mk - 1)$ for some integers $l$ and $k$. The reason is that there is only a very limited choice of possible forms of a solution, namely the two choices given by equations (5) and (6), in contrast to five more choices if $n_1 \equiv 0 \bmod m$ is allowed, see [7]. On the other hand the condition i) is by no means sufficient, because an integer solution puts also special requirements on the zeroth order solution $x_0, y_0, z_0$, e. g. $x_0 / n_0 \in \mathbb{N}$. See $n \equiv 7 \bmod 19$ as an example. 19 is not a multiple of $m = 5$, and $19 = mk - 1$ with $k = 4$, so a polynomal solution is not excluded by condition i) of the theorem. However, when parameters $s, r$ are calculated from condition ii) one gets $s = 5 + 7t, r = 13 + 19t$ with $t \in \mathbb{Z}_+$, which violates condition iii) for all integer $t$. One may as well argue from $5/7 = 1/x_0 + 1/y_0 + 1/z_0$ and calculate all possible solutions of this equation to see that there is no integer $x_0$ with $4 = k = x_0 / n_0$.

Another point is: Instead of solving a <u>nonlinear</u> diophantine equation for $x_0, y_0, z_0$ to sort out the parameters which lead to integer coefficients in (5) or (6) one just needs to determine parameters $k$ and $l$ from condition i) and solve a <u>linear</u> diophantine equation for $r$ and $s$. Since $(n_0, n_1) = 1$ divides $kl$, a solution to the linear equation ii) always exists. Given a special solution of ii), say $s_0, r_0$, one gets all solutions via $s = s_0 + n_0 t$ and $r = r_0 + n_1 t$, $t \in \mathbb{Z}_+$ and then singles out those ones which fulfill condition iii). If there are none, there is no polynomial solution to equation (2).

**References**

[1] A. Aigner, Brüche aus Summen von Stammbrüchen, *J. Angew. Math.* 214/215 (1964), 174-179.

[2] L. Bernstein, Zur Lösung der diophantischen Gleichung m/n=1/x+1/y+1/z insbesondere im Fall m=4, *Journal für die Reine und Angewandte Mathematik* Vol. 211 (1962), 1–10.

[3] P. Erdös and R.L. Graham, Old and new problems and results in combinatorial number theory, *Monographies de L'Enseignement Mathematique de Geneve* 28 (1980) pp. 30-44.

[4] M.Gionfriddo & E.Guardo (2021): A short proof of Erdös–Straus conjecture for every $n \equiv 13$ mod 24, *Journal of Interdisciplinary Mathematics* Vol 24, Issue 8, doi: 10.1080/09720502.2021.1907018.

[5] Graham, R.L. (2013). Paul Erdős and Egyptian Fractions. In: Lovász, L., Ruzsa, I.Z., Sós, V.T. (eds) Erdős Centennial. Bolyai Society Mathematical Studies, vol 25. Springer, Berlin, Heidelberg. doi: 10.1007/978-3-642-39286-3_9

[6] A. Schinzel, On sums of three unit fractions with polynomial denominators, *Functiones et Approximatio* XXVIII (2000), 187 – 194.

[7] B.R. Schuh, On polynomial solutions of the diophantine equation $m/n=1/x+1/y+1/z$, Preprint. Doi: 10.13140/RG.2.2.21739.52007